\documentclass[12pt]{amsart}

\usepackage{latexsym, amssymb, amsmath,ulem,soul,esint}

\usepackage{cancel}
\usepackage{latexsym, amssymb, amsmath}
\usepackage{soul,esint}
\usepackage{amsfonts, graphicx}
\usepackage{graphicx,color}

\usepackage{graphicx,color}
\newcommand{\be}{\begin{equation}}
\newcommand{\ee}{\end{equation}}
\newcommand{\beq}{\begin{eqnarray}}
\newcommand{\eeq}{\end{eqnarray}}

\newtheorem{thm}{Theorem}[section]

\newtheorem{lma}{Lemma}[section]
\newtheorem{prop}{Proposition}[section]

\newtheorem{claim}{Claim}[section]

\theoremstyle{remark}
\newtheorem{rem}{Remark}[section]
\numberwithin{equation}{section}
\def\tr{\operatorname{tr}}

\def\be{\begin{equation}}
\def\ee{\end{equation}}
\def\bee{\begin{equation*}}
\def\eee{\end{equation*}}

\def\K{K\"ahler }

\def\KR{K\"ahler-Ricci }
\def\Ric{\text{\rm Ric}}

\def\tr{\operatorname{tr}}

\def\e{\varepsilon}
\def\a{{\alpha}}
\def\b{{\beta}}

\def\ddb{\sqrt{-1}\partial\bar\partial}

\begin{document}

\title[]{On the H\"older estimate of K\"ahler-Ricci flow}

 \author{Jianchun Chu}
\address[Jianchun Chu]{Department of Mathematics, Northwestern University, 2033 Sheridan Road, Evanston, IL 60208}
\email{jianchun@math.northwestern.edu}

 \author{Man-Chun Lee$^1$}
\address[Man-Chun Lee]{Department of Mathematics, Northwestern University, 2033 Sheridan Road, Evanston, IL 60208; Mathematics Institute, Zeeman Building,
University of Warwick, Coventry CV4 7AL}
\email{mclee@math.northwestern.edu, Man.C.Lee@warwick.ac.uk}
\thanks{$^1$Research partially supported by NSF grant DMS-1709894 and EPSRC grant number P/T019824/1.}

\renewcommand{\subjclassname}{
  \textup{2010} Mathematics Subject Classification}

\date{\today}

\begin{abstract}
In this work, we study the H\"older regularity of the K\"ahler-Ricci flow on compact K\"ahler manifolds with semi-ample canonical line bundle. By adapting the method in the work of Hein-Tosatti on collapsing Calabi-Yau metrics, we obtain a uniform spatial $C^\alpha$ estimate along the K\"ahler-Ricci flow as $t\rightarrow +\infty$. %This improves the regularity of the K\"ahler-Ricci flow in the earlier works.
\end{abstract}


\maketitle

\markboth{Jianchun Chu and Man-Chun Lee}{}

\section{Introduction}
In this work, we study the normalized \KR flow which is a family of \K metrics satisfying
\begin{equation}\label{NKRF-equ}
\left\{
\begin{array}{ll}
\partial_t\omega(t)=-\Ric(\omega(t))-\omega(t);\\[1mm]
\omega(0)=\omega_0
\end{array}
\right.
\end{equation}
on a compact K\"ahler manifold $X$ with semi-ample canonical line bundle $K_X$, where $\omega_{0}$ is the initial K\"ahler metric. In this case, $X$ admits an Iitaka fibration structure given by a holomorphic map $f : X \to \Sigma \subset \mathbb{CP}^N$ with possibly singular fibers and possibly singular base manifold $\Sigma$. Let $S \subset \Sigma$ be the union of the set of singular values of $f$ and the singular set of $\Sigma$. The regular fibers $f^{-1}(z)$, where $z \in \Sigma \backslash S$, are Calabi-Yau manifolds. The complex dimension of $\Sigma$ is the Kodaira dimension of $X$. We focus on the case when $0 < \dim \Sigma < \dim X$, and we let $\dim_{\mathbb{C}} \Sigma = m$ and $\dim_{\mathbb{C}} X = m+n$, so that the Calabi-Yau fibers have complex dimension $n$.

Under semi-ample assumption, the canonical line bundle is nef and hence the flow exists on $X\times[0,+\infty)$ by the works of \cite{Cao1985,TianZhang2006,Tsuji1988}. The K\"ahler-Ricci flow under semi-ample canonical line bundle has been extensively studied by various authors \cite{FongZhang2015,FongYSZhang2019, Gill2014,GrossTosattiZhang2019,Jian2018,JianShi2019,SongTian2007,SongTian2012,SongTian2016,SongTianZhang2019,TianZLZ18,TosattiZhang2015,TosattiWeinkoveYang2018,YZhang2017,YZhang2018,YSZhang2019,ZZhang2014}. In \cite{SongTian2007,SongTian2012}, Song-Tian proved that the flow converges to a generalized K\"ahler-Einstein metric in the sense of measure on the base manifold $\Sigma$ as $t\rightarrow +\infty$. The generalized K\"ahler-Einstein metric $\omega_\Sigma$ satisfies $\Ric(\omega_\Sigma) = -\omega_\Sigma + \omega_{\text{WP}}$, where $\omega_{\textup{WP}}$ is the Weil-Petersson form which measures the variation of complex structures of the fibers. It was conjectured that the regularity of convergence can be improved to $C^\infty_{\textup{\textup{loc}}}(f^{-1}(\Sigma\backslash S))$-convergence. This conjecture is still open in general, although many progresses have been made. For instances,  Tosatti-Weinkove-Yang proved in \cite{TosattiWeinkoveYang2018} the $C^0_{\textup{\textup{loc}}}(f^{-1}(\Sigma\backslash S))$-convergence of the metric to the generalized K\"ahler-Einstein metric on the base manifold $\Sigma$. In \cite{FongLee2020}, the second named author and Fong considered the case when the generic fibres are biholomorphic to each other and developed a sharp parabolic Schauder estimate on cylinder using the idea of Hein-Tosatti in \cite{HeinTosatti2020}, and thus confirmed the above conjecture in the locally product case. More recently, Jian and Song \cite{JianSong2021} considered the case when $m+n=3$ and proved that the Ricci curvature is uniformly locally bounded. For further discussions, we refer interested readers to  \cite{Broder2020,FongLee2020,FongZhang2015,Gill2014,GrossTosattiZhang2019,HeinTosatti2015,JianSong2021,SongTianZhang2019,TianZLZ18,TosattiWeinkoveYang2018} and the references therein.

The main goal of this article is to establish the $C^\a_{\textup{\textup{loc}}}(f^{-1}(\Sigma\backslash S))$-convergence of the flow in the general setting where the fibres are not necessarily biholomorphic to each other.

\begin{thm}\label{intro-Holder-KRF}
Suppose that $(X,\omega_X)$ is a compact \K manifold with semi-ample canonical line bundle and $\omega(t)$ is a normalized \KR flow on $X$ defined by \eqref{NKRF-equ}. Then for any $\a\in (0,1)$ and any compact set $K\Subset X\setminus f^{-1}\left(S\right)$, there is $C(\a,K)$ such that for all $t\in [0,+\infty)$,
\begin{equation}
\|\omega(t)\|_{C^\a(K,\omega_X)}\leq C.
\end{equation}
\end{thm}

The proof of Theorem \ref{intro-Holder-KRF} adapts the idea in Hein-Tosatti's work \cite{HeinTosatti2020} on collapsing Calabi-Yau metrics, and some estimates for \KR flow in earlier works by Fong-Zhang \cite{FongZhang2015}, Song-Tian \cite{SongTian2016} and Tosatti-Weinkove-Yang \cite{TosattiWeinkoveYang2018}.

\vskip0.3cm

Acknowledgement: The authors are grateful to Valentino Tosatti and Hans-Joachim Hein for continuous support throughout the work. The authors would also like to thank Frederick Fong for useful discussion.

\section{Preliminary}\label{main-ingredient}
In this section, we collect some known results which will be used in deriving the $C^\a$ estimate.
\subsection{Local estimates for \KR flow}
It is known that when the \KR flow $\omega(t)$ is uniformly equivalent to a fixed \K metric, then $\omega(t)$ is bounded in $C^k_{\mathrm{loc}}$ for all $k\in \mathbb{N}$. When the \KR flow is with respect to a mildly varying family of complex structures, we need the following parabolic regularization property, which is a slight modification of the elliptic case \cite[Proposition 2.3]{HeinTosatti2020}.

\begin{prop}\label{KRF-localestimate}
For all $n,k,l_0\in \mathbb{N}$, $\a\in (0,1)$ and $A>1$, there are $\kappa(n,\a)$ and $C(k,n,\a,A,l_0)>0$ such that the following holds. Let $B_1(0)$ be the unit ball in $\mathbb{C}^n$ with the standard Euclidean metric $\omega_{\mathbb{C}^n}$ and $J$ be a complex structure on $B_1(0)$ such that
\begin{equation}
\|J-J_{\mathbb{C}^n}\|_{C^{1,\a}(B_1(0))}<\kappa, \quad \|J-J_{\mathbb{C}^n}\|_{C^{k,\a}(B_1(0))}\leq A.
\end{equation}
If $\omega(t)$ is a family of $J$-\K metrics satisfying
\begin{equation}
\partial_t \omega(t)=-\Ric(\omega(t))-l\omega(t)\;\;\text{on}\;\;Q_1(0)=B_1(0)\times [-1,0]
\end{equation}
for some $|l|\leq l_0$ and
\begin{equation}
 A^{-1}\omega_{\mathbb{C}^n}\leq \omega(t)\leq A\omega_{\mathbb{C}^n}\;\; \text{on}\;\;Q_1(0),
\end{equation}
then we have
\begin{equation}
\|\omega(t)\|_{C^{k,\a}(Q_{1/2}(0))}\leq C(k,n,\a,A,l_0).
\end{equation}
\end{prop}
\begin{proof}
By re-parametrization of time, we may assume $l_0=0$. Thanks to \cite[Proposition 2.2]{HeinTosatti2020}, we can find $J$-holomorphic coordinates on $B_{3/4}(0)$ which are close to the standard Euclidean coordinates in $C^{2,\a}$ and differing from them by a bounded amount in $C^{k+1,\a}$. In these $J$-holomorphic coordinates, applying the local estimates of \KR flow in \cite{ShermanWeinkove2012}, we obtain the required estimate.
\end{proof}

\subsection{Liouville theorems for Ricci-flat metric}
The following Liouville theorems for Ricci-flat metrics will play an important role in analyzing the blow-up model of the \KR flow.
\begin{thm}[\cite{RiebesehlFriedmar1984}]\label{LV-Euc}
Suppose that $\omega$ is a Ricci-flat \K metric on $\mathbb{C}^n$ such that
\begin{equation}
A^{-1}\omega_{\mathbb{C}^n} \leq \omega\leq A\omega_{\mathbb{C}^n} \;\;\text{on}\;\; \mathbb{C}^n
\end{equation}
for some $A>1$, then $\omega$ is constant.
\end{thm}

When the underlying manifold is a product of $\mathbb{C}^n$ and a compact Calabi-Yau \K manifold $Y$, the following Liouville theorem was proved by Hein \cite{Hein2019}, see also \cite{LiLiXi2020} for an alternative proof using mean value inequality.
\begin{thm}[\cite{Hein2019,LiLiXi2020}]\label{LV-prod}
Let $Y$ be a compact Calabi-Yau \K manifold with Ricci-flat metric $\omega_{Y}$. Suppose that $\omega$ is a Ricci-flat \K metric on $\mathbb{C}^n\times Y$ such that
\begin{equation}
A^{-1} (\omega_{\mathbb{C}^n}+\omega_Y) \leq \omega\leq A (\omega_{\mathbb{C}^n}+\omega_Y)\;\;\text{on}\;\; \mathbb{C}^n\times Y
\end{equation}
for some $A>1$ and $\omega$ is $d$-cohomologous to $\omega_{\mathbb{C}^n}+\omega_Y$, then $\omega$ is parallel with respect to $\omega_{\mathbb{C}^n}+\omega_Y$.
\end{thm}

\section{$C^\a$ estimate of \KR flow}
First, let us recall the setting: $X$ is a compact \K manifolds with semi-ample canonical line bundle, $f : X^{m+n} \to \Sigma^m \subset \mathbb{CP}^N$ is the corresponding Calabi-Yau fibration, and $S \subset \Sigma$ is the union of the set of singular values of $f$ and the singular set of $\Sigma$. By \cite{SongTian2012}, there exists a smooth \K metric $\omega_\Sigma$ on $\Sigma^n\setminus S$ satisfying the generalized K\"ahler-Einstein equation:
\begin{equation}\label{GKE-Sigma}
\Ric(\omega_\Sigma)=-\omega_\Sigma+\omega_{\mathrm{WP}},
\end{equation}
where $\omega_{\mathrm{WP}}$ is the smooth semi-positive Weil-Petersson form. By rescaling, we may assume $B_2(0) \Subset \Sigma \backslash S$, where $B_2(0)=B_{\mathbb{C}^m}(2)$ denotes the Euclidean ball. On $B_2(0)$, we have $\omega_{\mathbb{C}^m}=\omega_\Sigma+\ddb u$ for some $u\in C^{\infty}(B_2(0))$. For notational convenience, we still use $u$ and $\omega_{\Sigma}$ to denote their pull-backs to $B_{2}(0)\times Y$.

The map $f|_{f^{-1}(B_2(0))}:f^{-1}(B_2(0))\rightarrow B_2(0)$ is a proper surjective holomorphic map with $n$-dimensional Calabi-Yau fibers. For each $z\in B_2(0)$, we write $X_{z}=f^{-1}(z)$. For K\"ahler metric $\omega_0$ on $X$, using Yau's theorem \cite{Yau1978}, there is a unique \KR flat metric $\omega_{F,z}$ on each fibre $X_z$ which is cohomologous to $\omega_0|_{X_z}$. We may choose $\rho$ locally smoothly such that $\omega_F=\omega_0+\ddb \rho$ and $\omega_F|_{X_z}=\omega_{F,z}$. Following \cite{HeinTosatti2020}, we define the reference closed real $(1,1)$-form on $B_2(0)$ by
\begin{equation}
\omega^\natural(t)=(1-e^{-t})\omega_\infty+e^{-t}\omega_F,
\end{equation}
where $\omega_\infty=f^*\omega_{\mathbb{C}^m}$. Note that $\omega_{F}$ may be not positive definite along the base directions, and so $\omega^\natural(t)$ is not necessarily positive definite for all $t\geq0$. But the Cauchy-Schwarz inequality shows $\omega^\natural(t)$ is positive definite for $t$ sufficiently large, which is what we concern. For convenience, by translation of time, we always assume $\omega^\natural(t)>0$ for all $t\geq0$, and denote its associated Riemannian metric by $g^\natural(t)$.

Suppose that $\omega^\bullet(t)$ is the normalized \KR flow \eqref{NKRF-equ} defined on $X\times [0,+\infty)$. Let $g^\bullet(t)$ be the associated Riemannian metric. We use $\varphi$ to denote the solution of the following ODE:
\begin{equation}\label{MA}
\dot\varphi+\varphi=\log\frac{\omega^\bullet(t)^{m+n}}{e^{-nt}C^{m+n}_n\omega_F^n\wedge \omega_\Sigma^m}-u,\;\; \varphi(0)=\rho.
\end{equation}
Here we recall that $u$ is the function such that $\omega_{\infty}=\omega_\Sigma+\ddb u$. By taking $\ddb$ on both sides of \eqref{MA} and using $\ddb\log(\omega_F^n\wedge \omega_\Sigma^m)=\omega_{\Sigma}$ (see e.g., \cite[Section 5.7]{Tosatti18}), we see that
\begin{equation}
\omega^\bullet(t)=\omega^\natural(t)+\ddb\varphi(t)\;\;\text{on}\;\;\left(B_{2}(0)\times Y\right)\times[0,\infty).
\end{equation}
This reduces back to the setting analogous to its elliptic counterpart in  \cite{HeinTosatti2020}.

By Ehresmann's theorem and shrinking $B_2(0)$ slightly, $f$ is a smooth fiber bundle. We may choose a trivialization $\Phi:B_2(0)\times Y\rightarrow X$ where $Y=f^{-1}(0)$ such that $\Phi|_{\{0\}\times Y}:\{0\}\times Y\rightarrow f^{-1}(0)=Y$ is the identity map. On $B_{2}(0)\times Y$, we use $J^\natural$ to denote the complex structure inherited from $X$ via the trivialization $\Phi$. Then the projection $\mathrm{pr}_{\mathbb{C}^m}$ is a $J^\natural$-holomorphic submersion. Let $J_{Y,z}$ denote the restriction of $J^\natural$ to the $J^\natural$-holomorphic fibre $\{z\}\times Y$. Note that $\Phi^*\omega_{F,z}$ is a Ricci-flat $J_{Y,z}$-\K metric on $\{z\}\times Y$. Let $g_{Y,z}$ be its associated Riemannian metric on $\{z\}\times Y$. Extend it trivially to the product metric on $\mathbb{C}^m\times Y$ and define the product shrinking metric $g_{z}(t)=g_{\mathbb{C}^m}+e^{-t}g_{Y,z}$, which is \K with respect to $J_z=J_{\mathbb{C}^m}+J_{Y,z}$. By the trivialization $\Phi$, we may assume the above complex structures, metrics and the \KR flow are defined on $B_2(0)\times Y$. We will omit the trivialization for notational convenience.

The main objective of this section is to prove the following $C^\a$ estimate of the \KR flow.

\begin{thm}\label{Holder-KRF}
For all $\a\in (0,1)$, there is $C>0$ such that for all $t\in [0,+\infty)$,
\begin{equation}
\sup_{x=(z,y)\in B_{1/4}\times Y}\sup_{(x',t')\in B^{g_{z}(t)}((x,t),\frac{1}{8})}\frac{|\eta(x,t)-\mathbf{P}_{xx'}^{g_{z}(t)}(\eta(x',t))|_{g_z(t)}}{d^{g_{z}(t)}(x,x')^\a}\leq C,
\end{equation}
where $\eta=\ddb\varphi$ and $\mathbf{P}_{xx'}^{g^\bullet(t)}$ denotes the $g_{z}(t)$-parallel transport along the $g_z(t)$-geodesic from $x'$ to $x$.
\end{thm}

\begin{rem}\label{rem-shorttime}We note that the Ricci flow $g^\bullet(t)$ is uniformly equivalent to the reference shrinking Riemannian metric $g_z(t)$ by the work of \cite{FongZhang2015}. By standard parabolic regularity theory, the flow is $C^k$ regular in finite time. The major difficulties is the uniformity for all $t>0$. 
\end{rem}

To prove Theorem \ref{Holder-KRF}, we first need the following lower order estimates along the \KR flow.

\begin{lma}\label{lower-order-estimate}
Under the above setting, there are $C>0$ and $T>1$ such that for all $(x,t)\in B_{2}(0)\times [T,+\infty)$,
\begin{enumerate}\setlength{\itemsep}{1mm}
\item[(i)] $|R_{g^\bullet(t)}|\leq C$;
\item[(ii)] $C^{-1}\omega^\natural(t)\leq \omega^\bullet(t)\leq C\omega^\natural(t)$;
\item[(iii)]$|\nabla(\dot\varphi+\varphi+u)|_{g^\bullet(t)}\leq C$;
\item[(iv)] $ e^{nt}{\omega^\bullet(t)^{m+n}}\rightarrow {C^{m+n}_n \omega_F^n\wedge \omega_{\Sigma}^m }$ uniformly on $B_{2}(0)$ as $t\to+\infty$.
\end{enumerate}
\end{lma}
\begin{proof}
The uniform boundedness of scalar curvature follows directly from \cite[Theorem 1.1]{SongTian2016}. By \cite[Theorem 1.1]{FongZhang2015}, the flow $\omega^{\bullet}(t)$ is uniformly equivalent to the shrinking reference metric $(1-e^{-t})\omega_\Sigma+e^{-t}\omega_{0}$ and hence $\omega^\natural(t)$ on $B_{2}(0)$. These prove (i) and (ii).

By \cite[Proposition 3.1]{SongTian2016}, the function $|\nabla v|_{g^{\bullet}(t)}$ is uniformly bounded where $\displaystyle v=\log \frac{\omega^\bullet(t)^{m+n}}{e^{-nt}\Omega}$ and $\Omega$ is a smooth volume form on $X$. Since $\Omega$ and $\omega_F^n\wedge \omega_{\infty}^m$ are uniformly equivalent on $B_{2}(0)$, then (iii) follows. The volume convergence (iv) follows from \cite[Lemma 3.1]{TosattiWeinkoveYang2018}.
\end{proof}

Next, let us include some important observation from \cite{HeinTosatti2020}. Choosing a complex coordinate chart $(y^1,...,y^n)$ on $Y$, together with complex coordinate chart $(z^1,...,z^m)$ on $B_{2}(0)$, $(z,y)$ is a complex coordinate chart on $B_{2}(0)\times Y$. Since $\mathrm{pr}_{\mathbb{C}^m}$ is holomorphic with respect to $J^\natural$ and $J_{z_0}$, then
\begin{equation}\label{complexstructure 1}
(\mathrm{pr}_{\mathbb{C}^m})_{*}\circ(J^\natural-J_{z_0}) = J_{\mathbb{C}^{m}}\circ(\mathrm{pr}_{\mathbb{C}^m})_{*}-J_{\mathbb{C}^{m}}\circ(\mathrm{pr}_{\mathbb{C}^m})_{*} = 0.
\end{equation}
Thus, ignoring the distinction between these complex coordinates and their complex conjugates, \eqref{complexstructure 1} shows
\begin{equation}\label{complexstructure}
(J^\natural-J_{z_0})\Big|_{(z,y)}=A(z_0,z,y)dz\otimes \partial_y+B(z_0,z,y)dy\otimes \partial_y,
\end{equation}
where $A,B$ are smooth matrix-valued functions with $B(z_0, z_0, y) = 0$. Combining \eqref{complexstructure} with the definitions of $g_{z_0}(t)$ and $g^\natural(t)$, we see that
\begin{equation}\label{metricvariation}
\begin{split}
& \big(g^\natural(t)-g_{z_0}(t)\big)\Big|_{(z,y)} \\
= {} & e^{-t}\Big(C(z_0,z,y) dz\otimes dz+D(z_0,z,y)dz\otimes dy+E(z_0,z,y)dy\otimes dy\Big),
\end{split}
\end{equation}
where $C,D,E$ are smooth matrix valued functions with $E(z_0,z_0,y)=0$. Thanks to the factor $e^{-t}$ in \eqref{metricvariation} and the Cauchy-Schwarz inequality, we can find $C,T>1$ such that for all $(z,t)\in B\times [T,+\infty)$,
\begin{equation}\label{metrics equivalence}
C^{-1}g_{z}(t)\leq g^\natural(t)\leq Cg_z(t).
\end{equation}
By translation of time, we will assume $T=0$ since we only concern the behaviour of the \KR flow as $t\rightarrow +\infty$ by Remark~\ref{rem-shorttime}.

Now, we are in a position to prove Theorem \ref{Holder-KRF}.

\begin{proof}[Proof of Theorem~\ref{Holder-KRF}]
We will follow closely the argument in \cite{HeinTosatti2020} and adapt the argument in the parabolic setting. Let $B=B_1(0)\subset B_2(0)\Subset \Sigma\setminus S$. In the following, we say a constant is uniform if it is independent of $t$, and always use $C$ to denote a uniform constant.

For $x=(z,y)\in B\times Y$, we consider the function
\begin{equation}
\begin{split}\mu(x,t)&=\left(d^{g_{z}(t)}\left(x,\partial(B\times Y)\right)\right)^\a\\
&\quad \times \sup_{x'\in B^{g_{z}(t)}\left(x,\frac{1}{4}d^{g_{z}(t)}(x,\partial(B\times Y))\right)}\frac{|\eta(x,t)-\mathbf{P}^{g_{z}(t)}_{xx'}(\eta(x',t))|_{g_{z}(x,t)}}{d^{g_{z}(t)}(x,x')^\a}.
\end{split}
\end{equation}
where $\mathbf{P}^{g_{z}(t)}_{xx'}$ denotes the $g_{z}(t)$-parallel transport along the $g_z(t)$-geodesic from $x'$ to $x$. To prove Theorem \ref{Holder-KRF}, it suffices to show that 
\begin{equation}
\sup_{(B\times Y)\times [0,\infty)} \mu(x,t) \leq C
\end{equation}
for some $C>0$. Suppose on the contrary, we can find sequences $t_i\rightarrow +\infty$ and $x_i\in B\times Y$ such that
\begin{equation}\label{mu x i t i 1}
\mu(x_i,t_i)=\sup_{(B\times Y)\times [0,t_i]}\mu(x,t)\rightarrow +\infty.
\end{equation}
Let $x_i=(z_i,y_i)$ and $x_i'\in \overline{B^{g_{z_i}(t_i)}(x_i,\frac14 d^{g_{z_i}(t_i)}(x_i,\partial(B\times Y)))}$ be the point realizing the supremum in the definition of $\mu$. Define $\lambda_{i}$ by
\begin{equation}\label{lambda i}
\lambda_{i}^{\alpha} = \frac{|\eta(x_i,t_i)-\mathbf{P}^{g_{z_i}(t_i)}_{x_ix_i'}(\eta(x_i',t_i))|_{g_{z_i}(x_i,t_i)}}{d^{g_{z_i}(t_i)}(x_i,x_i')^\a},
\end{equation}
and then
\begin{equation}\label{mu x i t i 2}
\mu(x_i,t_i)= d^{g_{z_i}(t_i)}\left(x_i,\partial(B\times Y) \right)^\a \cdot \lambda_i^\a.
\end{equation}
By passing to subsequence, we may assume that
\begin{equation}\label{convergence-basepoint}
x_i\rightarrow x_\infty=(z_\infty,y_\infty)\in \overline{B}\times Y.
\end{equation}
Since $g_{z_{i}}(t_{i})=g_{\mathbb{C}^m}+e^{-t_{i}}g_{Y,z_{i}}$ and $t_{i}\rightarrow+\infty$, the metric $g_{z_{i}}(t_{i})$ is shrinking, and so $ d^{g_{z}(t_i)}\left(x,\partial(B\times Y)\right)$ is uniformly bounded from above. Combining this with \eqref{mu x i t i 1} and \eqref{mu x i t i 2}, we see that
\begin{equation}
\lambda_i\rightarrow +\infty\quad\text{as}\;\;i\rightarrow +\infty.
\end{equation}
Define the diffeomorphism $\Psi_i:B_{\lambda_i}\times Y\rightarrow B\times Y$ by
\begin{equation}
\Psi_i(z,y)=(\lambda_i^{-1}z,y),
\end{equation}
and pull back the complex structures, metrics, \KR flow and points to $B_{\lambda_i}\times Y$ via $\Psi_i$:
\begin{equation}\label{hat notations}
\left\{
\begin{array}{ll}
\hat J_i&=\Psi_i^*J_{z_i};\\[1mm]
\hat J_i^\natural&=\Psi_i^*J^\natural;\\[1mm]
\hat g_{i}(t)&=\lambda_i^2\Psi_i^*g_{z_i}(t_i+\lambda_i^{-2}t);\\[1mm]
\hat\omega_i^\natural(t)&=\lambda_i^2\Psi_i^*\omega^\natural(t_i+\lambda_i^{-2}t);\\[1mm]
\hat\eta_i(t)&=\lambda_i^2\Psi_i^*\eta(t_i+\lambda_i^{-2}t);\\[1mm]
\hat\omega_i^\bullet(t)&=\hat\omega_i^\natural(t)+\hat\eta_i(t);\\[1mm]
\hat x_i&=\Psi_i^{-1}(x_i);\\[1mm]
\hat x_i'&=\Psi_i^{-1}(x_i').
\end{array}
\right.
\end{equation}
Then $\hat g_i(t)$ is a Ricci-flat $\hat J_i$-\K product metric, $\hat\omega_i^\natural(t)$ is a semi-Ricci-flat $\hat J_i^\natural$-\K metric and $\hat\omega_i^\bullet(t)$ is a \KR flow with respect to $\hat J_i^\natural$. Let $\hat g_i^\natural(t)$ and $\hat g_i^\bullet(t)$ be the Riemannian metrics associated to $\omega_i^\natural(t)$ and $\hat\omega_i^\bullet(t)$, then we have the following properties which follow from Lemma \ref{lower-order-estimate}, \eqref{metrics equivalence} and \eqref{hat notations}:
\begin{equation}\label{Metric-equ-1}
\left\{
\begin{array}{ll}
C^{-1}\hat g_i(t)\leq \hat g_i^\bullet(t)\leq C\hat g_i(t);\\[1.5mm]
C^{-1}\hat g_i(t)\leq g_i^\natural(t) \leq C\hat g_i(t);\\[1.2mm]
\hat g_i(t)=g_{\mathbb{C}^m}+\lambda_i^2 e^{-t_i-\lambda_i^{-2}t}g_{Y,z_i};\\[1mm]
\hat \omega_i^\natural(t)=(1-e^{-t_i-\lambda_i^{-2}t})\omega_{\mathbb{C}^m}+\lambda_i^2e^{-t_i-\lambda_i^{-2}t}\Psi_i^*\omega_F.
\end{array}
\right.
\end{equation}
Recalling the definition of $\lambda_{i}$ \eqref{lambda i} and using \eqref{hat notations},
\begin{equation}
\begin{split}
\lambda_i^\a&=\frac{|\eta(x_i,t_i)-\mathbf{P}^{g_{z_i}(t_i)}_{x_ix_i'}(\eta(x_i',t_i))|_{g_{z_i}(x_i,t_i)}}{d^{g_{z_i}(t_i)}(x_i,x_i')^\a}=\frac{|\hat \eta_i(\hat x_i,0)-\mathbf{P}^{\hat g_{i}(0)}_{\hat x_i\hat x_i'}(\hat \eta_i(\hat x_i',0))|_{\hat g_{i}(\hat x_i,0)}}{d^{\hat g_{i}(0)}(\hat x_i,\hat x_i') ^\a}\cdot \lambda_i^\a,
\end{split}
\end{equation}
which implies
\begin{equation}\label{semi-1}
\frac{|\hat \eta_i(\hat x_i,0)-\mathbf{P}^{\hat g_{i}(0)}_{\hat x_i\hat x_i'}(\hat \eta_i(\hat x_i',0))|_{\hat g_{i}(\hat x_i,0)}}{d^{\hat g_{i}(0)}(\hat x_i,\hat x_i')^\a}=1.
\end{equation}
By \eqref{Metric-equ-1} and $\hat\omega_i^\bullet(t)=\hat\omega_i^\natural(t)+\hat\eta_i(t)$, the numerator of \eqref{semi-1} is uniformly bounded and hence the distance between $\hat x_i$ and $\hat x_i'$ with respect to $\hat g_i(0)$ is uniformly bounded:
\begin{equation}\label{bounded-spacetime}
d^{\hat g_{i}(0)}(\hat x_i,\hat x_i')\leq C.
\end{equation}

On the other hand, by the definition of $\mu$ and \eqref{mu x i t i 1},
\begin{equation}\label{complete-plimit}
\begin{split}
\left(d^{\hat g_i(0)}(\hat x_i,\partial(B_{\lambda_i}\times Y))\right)^\a = \mu(x_i,t_i) \rightarrow +\infty.
\end{split}
\end{equation}
This implies that the pointed limit space of
\begin{equation}
\left( B_{\lambda_i}\times Y,\hat g_i(0),\hat x_i\right)
\end{equation}
will be complete. By \eqref{bounded-spacetime}, we may assume $\hat x_i=(\hat z_i,\hat y_i)=(0,\hat y_i)$ modulo translations in the $\mathbb{C}^m$ factor. Write $\delta_i=\lambda_ie^{-t_i/2}$, then
\begin{equation}
\hat g_i(0)=g_{\mathbb{C}^m}+\delta_{i}^{2}g_{Y,z_i}.
\end{equation}
By passing to a subsequence, we may assume $\delta_i\rightarrow \delta_\infty\in [0,+\infty]$. From the behaviour of $\hat g_i(0)$, there are three cases to be considered:
\begin{enumerate}\setlength{\itemsep}{1mm}
\item[\bf(a)] $\delta_\infty=+\infty$;
\item[\bf(b)] $\delta_\infty\in (0,+\infty)$;
\item[\bf (c)] $\delta_\infty=0$.
\end{enumerate}
Before splitting into different cases, thanks to \eqref{complexstructure} and \eqref{convergence-basepoint}, we always have the following convenience of complex structures:
\begin{equation}\label{complexstructure-convergence}
\hat J_i,\;\hat J_i^\natural \rightarrow J_{\mathbb{C}^m}+J_{Y,z_\infty}\;\;\text{in}\;\; C^\infty_{\mathrm{loc}}(\mathbb{C}^m\times Y).
\end{equation}
\vskip0.3cm

{\bf Case (a): $\delta_\infty=+\infty$.} In this case, $\left( B_{\lambda_i}\times Y,\hat g_i(0),\hat x_i\right)$ converges to $(\mathbb{C}^{m+n},g_{\mathbb{C}^{m+n}},0)$ in the $C^\infty$-Cheeger-Gromov sense. More precisely, let $(\hat y^1,...,\hat y^n)$ be a holomorphic chart of $Y$ centred at $y_\infty$ with respect to the complex structure $J_{Y,z_\infty}$. We may assume $\hat y_i\in B_{\mathbb{C}^n}(1)\Subset Y$ and $\hat y_i\rightarrow \hat y_\infty=0$.

Consider the diffeomorphism $\Lambda_i:B_{\lambda_i}\times B_{\delta_i}\rightarrow B_{\lambda_i}\times B_1$ given by
$$\Lambda_i(z,y)=(z,\delta_i^{-1}y).$$
Then \eqref{complexstructure-convergence} shows the convergence of the background complex structures:
\begin{equation}\label{Case1:cs}
\Lambda_i^* \hat J_i,\;\Lambda_i^*\hat J_i^\natural \rightarrow J_{\mathbb{C}^{m+n}}\;\;\text{in}\;\;C^\infty_{\mathrm{loc}}(\mathbb{C}^{m+n}),
\end{equation}
and \eqref{metricvariation}, \eqref{Metric-equ-1} show the convergence of the background metrics:
\begin{equation}\label{Case1:m}
\Lambda_i^* \hat g_i(t),\;\Lambda_i^*\hat g_i^\natural(t) \rightarrow g_{\mathbb{C}^{m+n}}\;\;\text{in}\;\;C^\infty_{\mathrm{loc}}(\mathbb{C}^{m+n}\times (-\infty,0]).
\end{equation}

On the other hand, we write $\tilde \omega_i^\bullet(t)=\Lambda_i^*\hat \omega_i^\bullet(t)$, and then $\tilde \omega_i^\bullet(t)$ solves the approximated \KR flow with respect to the complex structure $\Lambda_i^* \hat J_i^\natural$:
\begin{equation}\label{Case1:aKRF}
\partial_t \tilde \omega_i^\bullet (t)=-\Ric(\tilde \omega_i^\bullet (t))-\lambda_i^{-2} \tilde \omega_i^\bullet (t).
\end{equation}
By \eqref{Metric-equ-1} and \eqref{Case1:m}, we obtain
\begin{equation}\label{Case1:m 1}
C^{-1}\omega_{\mathbb{C}^{m+n}}\leq \tilde \omega_i^\bullet(t)\leq C\omega_{\mathbb{C}^{m+n}}
\end{equation}
on any compact subset of $\mathbb{C}^{m+n}\times (-\infty,0]$ if $i$ is sufficiently large relative to the compact set. By \eqref{Case1:cs}, \eqref{Case1:m} and \eqref{Case1:m 1}, we may apply Proposition~\ref{KRF-localestimate} to obtain $C^\infty_{\mathrm{loc}}(\mathbb{C}^{m+n}\times (-\infty,0])$ estimate of $\tilde \omega_i^\bullet(t)$. Hence, $\tilde \omega_i^\bullet(t)$ converges to $\tilde \omega_\infty^\bullet(t)$ in $C^\infty_{\mathrm{loc}}(\mathbb{C}^{m+n}\times (-\infty,0])$ which remains uniformly equivalent to the Euclidean metric $\omega_{\mathbb{C}^{m+n}}$ for all $t\leq 0$.

Now, we follow similar argument in \cite{FongLee2020} to reduce the discussion back to the elliptic case. Using \eqref{Case1:aKRF}, $\tilde \omega_\infty^\bullet(t)$ solves the \KR flow:
\begin{equation}
\partial_t \tilde \omega_\infty^\bullet (t)=-\Ric(\tilde \omega_\infty^\bullet (t)).
\end{equation}
By Lemma \ref{lower-order-estimate}, the scalar curvature of the original \KR flow is bounded, and hence the scalar curvature after parabolic rescaling converges to $0$, i.e.,
\begin{equation}
R(\tilde \omega_\infty^\bullet(t))\equiv 0\;\;\text{on}\;\; \mathbb{C}^{m+n}\times (-\infty,0].
\end{equation}
Recalling the evolution of scalar curvature along the K\"ahler Ricci flow:
\begin{equation}
(\partial_{t}-\Delta_{\tilde \omega_\infty^\bullet (t)})R(\tilde \omega_\infty^\bullet (t))=2|\Ric(\tilde \omega_\infty^\bullet (t))|^2,
\end{equation}
we conclude that $\tilde \omega_\infty^\bullet (t)$ is Ricci-flat, and then $\tilde \omega_\infty^\bullet(t)\equiv \tilde \omega_\infty^\bullet(0)$ for all $t\leq0$. Combining this with Theorem \ref{LV-Euc}, $\tilde\omega_\infty^\bullet(t)$ is constant on $\mathbb{C}^{m+n}\times (-\infty,0]$.

Define the pull-back geometric quantities by
\begin{equation}
\left\{
\begin{array}{ll}
\tilde g_i(t)&=\Lambda_i^*\hat g_i(t);\\[1mm]
\tilde \eta_i(t)&=\Lambda_i^*\hat\eta_i(t);\\[1mm]
\tilde x_i&=\Lambda_i^{-1}(\hat x_i);\\[1mm]
\tilde x_i'&=\Lambda_i^{-1}(\hat x_i').
\end{array}
\right.
\end{equation}
By \eqref{semi-1} and \eqref{bounded-spacetime}, we see that
\begin{equation}\label{semi-3}
\left\{
\begin{array}{ll}
d^{\tilde  g_i(0)}(\tilde x_i,\tilde x_i')\leq  C;\\\\
\displaystyle\frac{|\tilde \eta_i(\tilde x_i,0)-\mathbf{P}^{\tilde g_i(0)}_{\tilde x_i\tilde x_i'}(\tilde\eta_i(\tilde x_i',0))|_{\tilde g_i(\tilde x_i,0)}}
{d^{\tilde g_i(0)}(\tilde x_i,\tilde x_i')^\a}=1.
\end{array}
\right.
\end{equation}
Since $d^{\tilde  g_i(0)}(\tilde x_i,\tilde x_i')$ is uniformly bounded, then there is a compact set $\Omega$ containing $\tilde x_i$ and $\tilde x_i'$. Note that $\Omega$ is independent of $i$. Using the convergence of reference complex structures \eqref{Case1:cs} and metrics \eqref{Case1:m}, the $C_{\mathrm{loc}}^{1}$ estimate of $\tilde \omega_i^\bullet(t)$ from Proposition \ref{KRF-localestimate} implies
\begin{equation}
|\tilde \eta_i(\tilde x_i,0)-\mathbf{P}^{\tilde g_i(0)}_{\tilde x_i\tilde x_i'}(\tilde\eta_i(\tilde x_i',0))|_{\tilde g_i(\tilde x_i,0)}
\leq Cd^{\tilde g_i(0)}(\tilde x_i,\tilde x_i').
\end{equation}
Combining this with \eqref{semi-3},
\begin{equation}\label{semi-9}
C^{-1}\leq d^{\tilde g_i(0)}(\tilde x_i,\tilde x_i')\leq  C.
\end{equation}
However, $\tilde \omega_i^\bullet(t)$ and $\Lambda_i^*\hat\omega_i^\natural(t)$ both converge to a constant real $(1,1)$-forms in $C^\infty_{\mathrm{loc}}(\mathbb{C}^{m+n}\times (-\infty,0])$. Using \cite[Remark 3.7]{HeinTosatti2020} and $\tilde \eta_i(t)=\tilde \omega_i^\bullet(t)-\Lambda_i^*\hat \omega_i^\natural(t)$, \eqref{semi-3} shows $d^{\tilde g_i(0)}(\tilde x_i,\tilde x_i')\rightarrow0$, which contradicts with \eqref{semi-9}.

\vskip0.3cm
{\bf Case (b): $\delta_\infty\in (0,+\infty)$.} In this case, the blowup limit is $\mathbb{C}^m\times Y$. We may assume $\delta_\infty=1$ without loss of generality. This case is similar to {\bf Case (a)} except that we don't need to apply an additional diffeomorphism. Indeed, similar to \eqref{complexstructure-convergence}, \eqref{metricvariation} shows
\begin{equation}\label{case2:convergences}
\hat g_i(t),\;\hat g_i^\natural(t)\rightarrow g_{\mathbb{C}^m}+g_{Y,z_\infty} \quad\text{in}\;\; C^\infty_{\mathrm{loc}}(\mathbb{C}^m\times Y\times (-\infty,0]).
\end{equation}
Thanks to \eqref{Metric-equ-1}, \eqref{complexstructure-convergence} and \eqref{case2:convergences}, we may apply Proposition~\ref{KRF-localestimate} to obtain $C_{\mathrm{loc}}^{\infty}(\mathbb{C}^m\times Y\times (-\infty,0])$ estimate of $\hat\omega_i^\bullet(t)$. Hence, $\hat\omega_i^\bullet(t)$ converges to $\hat\omega_\infty^\bullet(t)$ in $C_{\mathrm{loc}}^{\infty}(\mathbb{C}^m\times Y\times (-\infty,0])$. Then for each $t\leq0$, $\hat\omega_\infty^\bullet(t)$ is cohomologous and uniformly equivalent to the product \K form $\omega_{\mathbb{C}^m}+\omega_{Y,z_\infty}$, and is \K with respect to the product complex structure $J_{\mathbb{C}^m}+J_{Y,z_\infty}$. As in {\bf Case (a)}, since the scalar curvature before applying diffeomorphism $\Psi_{i}$ is uniformly bounded, the same argument showed that $\hat \omega^\bullet_\infty(t)=\hat \omega^\bullet_\infty(0)$ is Ricci-flat. Thanks to the uniform equivalence of metrics \eqref{Metric-equ-1} and Theorem \ref{LV-prod}, $\hat \omega^\bullet_\infty(t)$ is parallel to $\omega_{\mathbb{C}^m}+\omega_{Y,z_\infty}$. But this will contradict with \eqref{semi-1} by the same argument of {\bf Case (a)}.

\vskip0.3cm
{\bf Case (c): $\delta_\infty=0$.} In this case, the blowup limit is $\mathbb{C}^m$. For $\hat x=(\hat z, \hat y)\in B_{\lambda_i}\times Y$ and $z=\lambda_i^{-1}\hat z$, we write
\begin{equation}
(\hat g_i)_{\hat z}(t)=\lambda_i^2\Psi_i^*g_{z}(t_i+\lambda_i^{-2}t)=g_{\mathbb{C}^m}+\lambda_i^2e^{-t_i-\lambda_i^{-2}t}g_{Y,z}.
\end{equation}
By the definition of $\mu$ and \eqref{semi-1}, for all $\hat x=(\hat z, \hat y)\in B_{\lambda_i}\times Y$ and
$\hat x'\in B^{(\hat g_i)_{\hat z}(0)}(\hat x,\frac{1}{4}d^{(\hat g_i)_{\hat z}(0)}(\hat x,\partial(B_{\lambda_i}\times Y)))$, we have
\begin{equation}
\begin{split}
& \quad\left(d^{(\hat g_i)_{\hat z}(0)}(\hat x,\partial(B_{\lambda_i}\times Y))\right)^{\alpha}
\frac{|\hat \eta_i(\hat x,0)-\mathbf{P}^{(\hat g_i)_{\hat z}(0)}_{\hat x\hat x'}(\hat\eta_i(\hat x',0))|_{(\hat g_i)_{\hat z}(\hat x,0)}}
{d^{(\hat g_i)_{\hat z}(0)}(\hat x,\hat x')^\a} \\[1mm]
& \leq \mu(x_i,t_i) = \left(d^{\hat g_i(0)}(\hat x_i,\partial(B_{\lambda_i}\times Y))\right)^\a,
\end{split}
\end{equation}
which implies
\begin{equation}
\frac{|\hat \eta_i(\hat x,0)-\mathbf{P}^{(\hat g_i)_{\hat z}(0)}_{\hat x\hat x'}(\hat\eta_i(\hat x',0))|_{(\hat g_i)_{\hat z}(\hat x,0)}}
{d^{(\hat g_i)_{\hat z}(0)}(\hat x,\hat x')^\a}
\leq \left(\frac{d^{\hat g_i(0)}(\hat x_i,\partial(B_{\lambda_i}\times Y))}{d^{(\hat g_i)_{\hat z}(0)}(\hat x,\partial(B_{\lambda_i}\times Y))}\right)^\a.
\end{equation}
Recalling that metrics $\hat g_i(0)$ and $(\hat g_i)_{\hat z}(0)$ are uniformly equivalent,
\begin{equation}\label{semi-6}
\frac{|\hat \eta_i(\hat x,0)-\mathbf{P}^{(\hat g_i)_{\hat z}(0)}_{\hat x\hat x'}(\hat\eta_i(\hat x',0))|_{(\hat g_i)_{\hat z}(\hat x,0)}}
{d^{(\hat g_i)_{\hat z}(0)}(\hat x,\hat x')^\a}
\leq C\left(\frac{d^{\hat g_i(0)}(\hat x_i,\partial(B_{\lambda_i}\times Y))}{d^{\hat g_i(0)}(\hat x,\partial(B_{\lambda_i}\times Y))}\right)^\a.
\end{equation}
For any $R>1$, $\hat x\in B^{\hat g_i(0)}(\hat x_i,R)$ and sufficiently large $i$, \eqref{complete-plimit} and \eqref{semi-6} show
\begin{equation}\label{semi-2}
\begin{split}
& \frac{|\hat \eta_i(\hat x,0)-\mathbf{P}^{(\hat g_i)_{\hat z}(0)}_{\hat x\hat x'}(\hat\eta_i(\hat x',0))|_{(\hat g_i)_{\hat z}(\hat x,0)}}
{d^{(\hat g_i)_{\hat z}(0)}(\hat x,\hat x')^\a} \\
\leq {} & C\left(\frac{d^{\hat g_i(0)}(\hat x_i,\partial(B_{\lambda_i}\times Y))}
{d^{\hat g_i(0)}(\hat x_i,\partial(B_{\lambda_i}\times Y))-d^{\hat g_i(0)}(\hat x,\hat x_i)}\right)^\a
\leq C_{R}.
\end{split}
\end{equation}
Combining this with \cite[Lemma 3.6]{HeinTosatti2020}, $\hat \omega_i^\bullet(0)$ have uniformly bounded $C^{\a}$ norm with respect to any fixed non-collapsing reference metric when $i$ is sufficiently large. Thus, we may assume $\hat\omega_i^\bullet(0)\rightarrow \hat\omega_\infty^\bullet(0)$ in $C^{\b}_{\mathrm{loc}}(\mathbb{C}^m\times Y)$ for any $\b<\a$, where $\hat\omega_\infty^\bullet(0)\in C^{\a}_{\mathrm{loc}}(\mathbb{C}^m\times Y)$ which satisfies the following:
\begin{enumerate}\setlength{\itemsep}{1mm}
\item[(i)] $\hat\omega_\infty^\bullet(0)$ is a section of $\mathrm{pr}^*_{\mathbb{C}^m}(\Lambda^{1,1}\mathbb{C}^m)$ uniformly equivalent to $\omega_{\mathbb{C}^m}$;
\item[(ii)] $\hat\omega_\infty^\bullet(0)$ is $g_{Y,z_\infty}$-parallel in the fiber directions;
\item[(iii)] $\hat\omega_\infty^\bullet(0)$ is weakly closed.
\end{enumerate}
The conclusions (i) and (ii) follows from \eqref{Metric-equ-1} and \eqref{semi-2}. The conclusion (iii) is clear due to the uniform convergence and the fact that $\hat\omega_i^\bullet(0)$ is closed. Hence, $\hat\omega_\infty^\bullet(0)$ is the pull-back under $\mathrm{pr}^*_{\mathbb{C}^m}$ of a weakly closed $(1,1)$-form of class $C^{\a}_{\mathrm{loc}}$ on $\mathbb{C}^m$.

To derive contradiction in {\bf Case (c)}, we will prove three claims. First, we rule out $d^{\hat g_i(0)}(\hat x_i,\hat x_i')\rightarrow 0$.
\begin{claim}\label{case3:noncoll}
There is $\e>0$ such that for all $i\in \mathbb{N}$,
\begin{equation}
d^{\hat g_i(0)}(\hat x_i,\hat x_i')\geq \e.
\end{equation}
\end{claim}

\begin{proof}[Proof of Claim~\ref{case3:noncoll}]
Suppose on the contrary, by passing to subsequence, we may assume
\begin{equation}
d_i := d^{\hat g_i(0)}(\hat x_i,\hat x_i')\rightarrow 0.
\end{equation}
Consider the diffeomorphism $\Lambda_i:B_{d_i^{-1}\lambda_i}\times Y\rightarrow B_{\lambda_i}\times Y$ given by
\begin{equation}
\Lambda_i(z,y)=(d_i z,y),
\end{equation}
and define the pull-back geometric quantities by
\begin{equation}
\left\{
\begin{array}{ll}
\tilde J_i&=\Lambda_i^*\hat J_{i};\\[1mm]
\tilde J_i^\natural&=\Lambda_i^*\hat J_i^\natural;\\[1mm]
\tilde g_{i}(t)&=d_i^{-2}\Lambda_i^*\hat g_{i}(d_i^2t);\\[1mm]
\tilde\omega_i^\natural(t)&=d_i^{-2}\Lambda_i^*\hat \omega_i^\natural(d_i^{2}t);\\[1mm]
\tilde\eta_i(t)&=d_i^2\Lambda_i^*\hat \eta_i(d_i^{-2}t);\\[1mm]
\tilde\omega_i^\bullet(t)&=\tilde\omega_i^\natural(t)+\tilde\eta_i(t);\\[1mm]
\tilde x_i&=\Lambda_i^{-1}(\hat x_i);\\[1mm]
\tilde x_i'&=\Lambda_i^{-1}(\hat x_i').
\end{array}
\right.
\end{equation}
Then
\begin{equation}\label{normal-dist}
d^{\tilde g_i(0)}(\tilde x_i,\tilde x_i')=1
\end{equation}
and
\begin{equation}\label{ref-1}
\tilde\omega_i^\natural(t)
=(1-e^{-s})\omega_{\mathbb{C}^m}+e^{-s}d_i^{-2}\lambda_i^2\Lambda_i^*\Psi_i^* \omega_F,
\end{equation}
where $s=t_i+d_i^2\lambda_i^{-2}t$. Hence, we can rewrite the equation \eqref{MA} as
\begin{equation}\label{1-MA}
(\tilde\omega_i^\natural+\tilde \eta_i)^{m+n}=e^{\tilde G_i+\tilde H_i}(\tilde\omega_i^\natural)^{m+n},
\end{equation}
where
\begin{equation}
\left\{
\begin{array}{ll}
\tilde G_i(t)&=\Lambda_i^* \Psi_i^* (\dot\varphi+\varphi+u)(s);\\
\tilde H_i(t)&=\displaystyle\log\frac{C^{m+n}_n\omega_{\Sigma}^m\wedge(e^{-s}d_i^{-2}\lambda_i^2\Lambda_i^*\Psi_i^* \omega_F)^n}
{\left((1-e^{-s})\omega_{\mathbb{C}^m}+e^{-s} d_i^{-2}\lambda_i^2\Lambda_i^*\Psi_i^* \omega_F \right)^{m+n}}.
\end{array}
\right.
\end{equation}

Using \eqref{Metric-equ-1} and \eqref{semi-1},
\begin{equation}\label{semi-5}
\left\{
\begin{array}{ll}
|\tilde\eta_i(0)|_{\tilde g_i(0)}\leq C;\\\\
\displaystyle\frac{|\tilde \eta_i(\tilde x_i,0)-\mathbf{P}^{\tilde g_i(0)}_{\tilde x_i\tilde x_i'}(\tilde\eta_i(\tilde x_i',0))|_{\tilde g_i(\tilde x_i,0)}}
{d^{\tilde g_i(0)}(\tilde x_i,\tilde x_i')^\a}=d_i^\a.
\end{array}
\right.
\end{equation}
For $\tilde x=(\tilde z,\tilde y)\in B_{d_i^{-1}\lambda_i}\times Y$ and $z=d_i\lambda_i^{-1}\tilde z$, we write
\begin{equation}
(\tilde g_i)_{\tilde z}(t)=d_i^{-2}\lambda_i^2\Lambda_i^*\Psi_i^*g_{z}(t_i+d_i^2\lambda_i^{-2}t).
\end{equation}
Using \eqref{semi-2}, after pulling back via the diffeomorphism $\Lambda_i$, we conclude that if $i$ sufficiently large,
\begin{equation}\label{semi-4}
\sup_{\tilde x ,\tilde x' \in B^{\tilde g_i(0)}\left(\tilde x_i,d_i^{-1}\right)}
\frac{|\tilde\eta_i(\tilde x,0)-\mathbf{P}^{(\tilde g_i)_{\tilde z}(0)}_{\tilde x\tilde x'}(\tilde \eta_i(\tilde x',0))|_{(\tilde g_i)_{\tilde z}(\tilde x,0)}}
{d^{\tilde g_i(0)}(\tilde x,\tilde x')^\a} \leq Cd_i^\a.
\end{equation}
Decompose $\tilde \eta_i(0)=\tilde\eta_i^\natural+\tilde\eta_i'$ where $\tilde\eta_i^\natural$ is the unique $\tilde g_i(0)$-parallel $(1,1)$-form pulled back from $\mathbb{C}^m$ such that $\tilde \eta_i^\natural(\tilde x_i)$ is the $\tilde g_i(\tilde x_i,0)$ orthogonal projection of $\tilde \eta_i(\tilde x_i,0)$ onto $\mathrm{pr}^*_{\mathbb{C}^m}(\Lambda^{1,1}\mathbb{C}^m)|_{\tilde x_i}$. Applying the proof of \cite[(5.36)]{HeinTosatti2020} by freezing the time at $0$ and using \eqref{semi-4}, we can find a constant $C>0$ such that for sufficiently large $i$,
\begin{equation}\label{semi-7}
|d_i^{-\a}\tilde\eta_i'(\tilde x_i)|_{\tilde g_i(\tilde x_i,0)}\leq C (d_i^{-1}\delta_i)^\a.
\end{equation}

Write $\e_i=d_i^{-1}\delta_i$. By passing to a subsequence, we may assume $\e_i\rightarrow \e_\infty\in [0,+\infty]$. By considering the behaviour of $\tilde g_i^\natural(0)=g_{\mathbb{C}^{m}}+\e_i^2 g_{Y,z_i}$, there are three distinct cases to be considered: $\e_\infty=+\infty$; $\e_\infty\in(0,+\infty)$; $\e_\infty=0$. The motivation of the above discussion is to pass $d_i^{-\a}\tilde\eta_i'$ to a limiting $\ddb $ exact $(1,1)$-form on $\mathbb{C}^{m+n},\mathbb{C}^m\times Y$ or $\mathbb{C}^m$ that is $O(r^\a)$ as $r\rightarrow+\infty$ and not parallel, where $r$ denotes the corresponding distance function in different settings.

Before analyzing these three cases, we collect some useful estimates which hold in all cases. Using \eqref{semi-4}, \eqref{semi-7} and $d_i\rightarrow 0$, for each $R>0$ and $i$ sufficiently large, we have
\begin{equation}\label{jet-1}
\sup_{\tilde x\in B^{\tilde g_i(0)}(\tilde x_i,R)}|\tilde\eta_i'(\tilde x)|_{\tilde g_i(\tilde x,0)}\leq C\delta_i^\a+Cd_i^\a R^\a.
\end{equation}
Following the argument of \cite[(5.41)]{HeinTosatti2020}, we define $\tilde \omega_i^\#=\tilde \omega_i^\bullet(0)-\tilde\eta_i'$. Since $\delta_i,d_i\rightarrow0$, then $\tilde \omega_i^\#$ is \K on $B^{\tilde g_i(0)}(\tilde x_i,R)$ for sufficiently large $i$. Moreover, its associated Riemannian metric $\tilde g_i^\#$ is uniformly equivalent to $\tilde g_i(0)$. At time $0$, we expand the Monge-Amp\`ere equation \eqref{1-MA} as
\begin{equation}\label{2-MA}
\tr_{\tilde \omega_i^\#} \tilde\eta_i'+\sum_{i=2}^{m+n}C^{m+n}_i\frac{(\tilde\eta_i')^i\wedge (\tilde\omega_i^\#)^{m+n-i}}{(\omega_i^\#)^{m+n}}
=e^{\tilde G_i(0)+\tilde H_i(0)} \frac{(\tilde\omega_i^\natural(0))^{m+n}}{(\tilde\omega_i^\#)^{m+n}}-1.
\end{equation}
Write $e^{\tilde K_i}-1$ for the right hand side.

\bigskip

{\bf Subclaim:} There is $C>0$ such that the following holds. For all $R>0$, there is $N$ such that for all $i>N$ and $\tilde x'\in B^{\tilde g_i(0)}(\tilde x_i,R)$,
\begin{enumerate}\setlength{\itemsep}{1mm}
\item[(i)] for $j\geq2$, $d_i^{-\a} |(\tilde\eta_i'(\tilde x_i))^j-\mathbf{P}^{\tilde g_i(0)}_{\tilde x_i\tilde x'}((\tilde \eta_i'(\tilde x'))^j)|_{\tilde g_i(\tilde x_i,0)} \leq C(\delta_i^\a +d_i^\a R^\a)R^\a$;
\item[(ii)] $d_i^{-\a}|\tilde\omega_i^\#(\tilde x_i)-\mathbf{P}^{\tilde g_i(0)}_{\tilde x_i\tilde x'}(\tilde \omega_i^\#(\tilde x'))|_{\tilde g_i(\tilde x_i,0)}\leq Cd_i^{1-\a} \lambda_i^{-1}R$;
\item[(iii)] $d_i^{-\a}|e^{\tilde K_i(\tilde x_i,0)}-e^{\tilde K_i(\tilde x',0)}|\leq Cd_i^{1-\a}\lambda_i^{-1} R$.
\end{enumerate}
\begin{proof}[Proof of Subclaim]
The proofs of (i) and (ii) are identical to that of \cite[(5.42), (5.43)]{HeinTosatti2020} by freezing the time at $0$ and using \eqref{complexstructure}, \eqref{metricvariation}, \eqref{semi-4} and \eqref{jet-1}.

For (iii), the contributions from $\tilde H_i(0)$ and ${(\tilde\omega_i^\natural(0))^{m+n}}/{(\tilde\omega_i^\#)^{m+n}}$ are done in the proof of \cite[(5.44)]{HeinTosatti2020}. It suffices to handle the term $\tilde G_i(0)$. We denote $v=\dot\varphi(0)+\varphi(0)+u$. For $\tilde x'$ such that $d_{\tilde g_i(0)}(\tilde x_i,\tilde x')<R$ and sufficiently large $i$,
\begin{equation}
\begin{split}
 |\tilde G_i(\tilde x_i,0)-\tilde G_i(\tilde x',0)| 
= {} & |v\left(\Psi_i^*\circ\Lambda_i(\tilde x_i) \right)-v\left(\Psi_i^*\circ\Lambda_i(\tilde x') \right)|\\
= {} & |v( z_i , y_i)-v(\lambda_i^{-1} d_i \tilde z', \tilde y')|\\
\leq {} & |\nabla v|_{L^\infty(B^{g(t_i)}(x_i,1))}\cdot d^{g(t_i)}\left((z_i,y_i),(\lambda_i^{-1}d_i \tilde z',\tilde y') \right).
\end{split}
\end{equation}
Combining this with Lemma \ref{lower-order-estimate}, \eqref{metricvariation} and \eqref{Metric-equ-1}, we obtain
\begin{equation}
|\tilde G_i(\tilde x_i,0)-\tilde G_i(\tilde x',0)|
\leq Cd^{g(t_i)}\left((z_i,y_i),(\lambda_i^{-1}d_i \tilde z',\tilde y') \right)
\leq Cd_i \lambda_i^{-1} R.
\end{equation}
Then the required estimate follows since $\tilde G_i$ is bounded by Lemma \ref{lower-order-estimate}.
\end{proof}

As a consequence of the Subclaim, we conclude that there is $C>0$ such that for $R>0$, $\tilde x'\in B^{\tilde g_i(0)}(\tilde x_i,R)$ and sufficiently large $i$,
\begin{equation}\label{semi-8}
d_i^{-\a}|\tr_{\tilde\omega_i^\#}\tilde\eta_i'(\tilde x_i)-\tr_{\tilde\omega_i^\#}\tilde\eta_i'(\tilde x')|
\leq C(\delta_i^\a+d_i^{\a}R^\a)R^\a+ Cd_i^{1-\a}\lambda_i^{-1}R,
\end{equation}
which is analogous to \cite[(5.48)]{HeinTosatti2020}. Applying the argument of \cite[5.3.1 Claim 1]{HeinTosatti2020} to derive contradiction in different cases: $\e_i\rightarrow +\infty$; $\e_i\rightarrow \e_\infty\in (0,+\infty)$; $\e_i\rightarrow 0$, we obtain Claim \ref{case3:noncoll}.
\end{proof}

\begin{claim}\label{case3:nonparallel}
There are two points $\hat z,\hat z'\in \mathbb{C}^m$ such that $\hat\omega_\infty^\bullet(\hat z,0)\neq \hat\omega_\infty^\bullet(\hat z',0)$.
\end{claim}
\begin{proof}
[Proof of Claim~\ref{case3:nonparallel}] The proof is identical to \cite[5.3.2 Claim 2]{HeinTosatti2020} by freezing time at $0$ and using Claim \ref{case3:noncoll}, \eqref{semi-1} and \eqref{bounded-spacetime}.
\end{proof}

\begin{claim}\label{case3:parallel}
The $C^\a$ \K current $\hat \omega_\infty^\bullet(0)$ on $\mathbb{C}^m$ is parallel with respect to Euclidean metric.
\end{claim}
\begin{proof}
[Proof of Claim~\ref{case3:parallel}]
The proof is standard and very similar to that in \cite{FongLee2020, HeinTosatti2020}. It suffices to show that
\begin{equation}\label{case3:volume-constant}
\hat\omega^\bullet_\infty(0)=c\omega_{\mathbb{C}^m}^m\;\;\text{weakly on}\;\;\mathbb{C}^{m}
\end{equation}
for some constant $c>0$. Indeed, since $\hat\omega_\infty^\bullet(0)\in C^{\a}_{\mathrm{loc}}(\mathbb{C}^m)$ is closed, then it has potential in $C^{2,\a}_{\mathrm{loc}}(\mathbb{C}^m)$. Combining this with \eqref{case3:volume-constant} and standard elliptic bootstrapping argument, $\hat\omega^\bullet_\infty(0)$ is smooth. Recalling that $\hat \omega_\infty^\bullet(0)$ is uniformly equivalent to $\omega_{\mathbb{C}^m}$, Theorem \ref{LV-Euc} implies $\hat\omega_\infty^\bullet(0)$ is constant, which complete the proof of Claim~\ref{case3:parallel}.

Next, we prove \eqref{case3:volume-constant}. Write $s=t_{i}+\lambda_{i}^{-2}t$ and $\psi_i(t)=\lambda_i^2 \Psi_i^* \varphi(s)$. Using \eqref{MA} and \eqref{hat notations}, we have
\begin{equation}\label{3-MA}
\left\{
\begin{array}{ll}
\hat\omega_i^\bullet(t)&=(1-e^{-s})\omega_{\mathbb{C}^m}+e^{-\lambda_i^{-2}t}\delta_i^2\Psi_i^* \omega_F +\ddb \psi_i(s);\\\\
\hat\omega_i^\bullet(t)^{m+n}&=\displaystyle C^{m+n}_n e^{\Psi_i^* (\dot\varphi(s)+\varphi(s)+u)}\delta_i^{2n} (\Psi_i^* \omega_F )^n\wedge \omega_{\mathbb{C}^m}^m.
\end{array}
\right.
\end{equation}
where $\partial$ and $\bar\partial$ are with respect to $\hat J_i^\natural$. Using $z_i\rightarrow z_\infty\in \overline{B}$ and \eqref{complexstructure}, we have the following convergence:
\begin{equation}\label{case3:convergences}
\hat J_i^\natural \rightarrow  J_{\mathbb{C}^m}+J_{Y,z_\infty},\;\; \Psi_i^* \omega_F\rightarrow  \omega_{Y,z_\infty}
\;\;\text{in}\;\; C^\infty_{\mathrm{loc}}(\mathbb{C}^m\times Y).
\end{equation}
Then \eqref{case3:volume-constant} follows from \eqref{3-MA}, \eqref{case3:convergences} and same argument of \cite{FongLee2020} which is based on \cite{HeinTosatti2020}.
\end{proof}

It is clear that Claim \ref{case3:nonparallel} contradicts with Claim \ref{case3:parallel}. This completes the proof of Case (c), and hence Theorem \ref{Holder-KRF}.
\end{proof}

\begin{proof}[Proof of Theorem~\ref{intro-Holder-KRF}]
The $C^0$ estimate of $\omega(t)$ follows from \cite[Theorem 1.1]{FongZhang2015}. Recalling $g_{z}(t)=g_{\mathbb{C}^m}+e^{-t}g_{Y,z}$, it is immediate that $\mathbf{P}^{g_z(t)}=\mathbf{P}^{g_z(0)}$ and $d^{g_z(t)}\leq d^{g_z(0)}$ and $|T|_{g_z(t)}\geq |T|_{g_z(0)}$ for all $t\geq0$ and contravariant tensor $T$. Therefore, Theorem \ref{Holder-KRF} implies a uniform bound on the $g_z(0)$-H\"older quotient of the \KR flow $\omega(t)$ for any $x,x'\in X$. Thanks to \cite[Lemma 3.6]{HeinTosatti2020}, we obtain the required H\"older estimate.
\end{proof}

\end{document}